\documentclass[12pt]{article}

\begin{document}

\begin{titlepage}
\title{\bf Applications over Complex Lagrangians}
\author{ Mehmet Tekkoyun \footnote{tekkoyun@pau.edu.tr} \\
 {\small Department of Mathematics, Pamukkale University,}\\
{\small 20070 Denizli, Turkey}\\
 Erdal  \"{O}z\"{u}saglam \footnote{materdalo@gmail.com} \\
 {\small Department of Mathematics, Aksaray University,}\\
{\small Aksaray, Turkey}\\
Ali G\"{o}rg\"{u}l\"{u} \footnote{agorgulu@ogu.edu.tr} \\
 {\small Department of Mathematics, Eski\c{s}ehir Osmangazi University,}\\
{\small Eski\c{s}ehir, Turkey}}
\date{\today}
\maketitle
\begin{abstract}
In this paper, Lagrangian formalisms of Classical Mechanics was
deduced on Kaehlerian manifold being geometric model of a
generalized Lagrange space.Then, it was given two applications of
complex Euler-Lagrange equations on mechanics system.\\
 {\bf Key
words:} Complex and Kaehlerian manifold, Lagrangian systems, Maple.
\end{abstract}
\end{titlepage}
\textbf{1.} \textbf{Introduction}\newline
As well known, modern differential geometry provides a suitable fields for
studying Lagrangian theory of Classical Mechanics. This is easily shown by
numerous articles and books \cite{Deleon, Norbury, Cabar, Gibbons, Tekkoyun,
Tekkoyun1} and there in. Therefore the dynamics of a Lagrangian system is
determined by a suitable vector field $X$ defined on the tangent bundle of a
given configuration space-manifold. If one takes an configuration manifold $M
$ and a regular Lagrangian function $L$\textbf{\ }on tangent bundle $TM$
then it is seen that there is an unique vector field $X$ on $TM$ such that
\begin{equation}
\begin{array}{l}
i_{X}\omega _{L}=dE_{L}
\end{array}
\label{1.1}
\end{equation}
where $\omega _{L}$ is the symplectic form and $E_{L}$ is energy associated
to $L$. The vector field $X$ is a second order differential
equation(semispray) since its integral curves are the solutions of the
Euler-Lagrangian equations. Here, we present the complex Euler-Lagrange
equations on Kaehlerian manifold being geometric model of a generalized
Lagrange space and to derive complex Euler-Lagrange equations on two
physical problems using Maple\cite{Celik}. Hereafter, all mappings and
manifolds are assumed to be differentiable of class $C^{\infty }$ and the
sum is taken over repeated indices$.$ Also, we denote by $\mathcal{F}(TM)$
the set of complex functions on $TM,$ by $\chi (TM)$ the set of complex
vector fields on $TM$ and by $\wedge ^{1}(TM)$ the set of complex 1-forms on
$TM.$ $1\leq i\leq n.$ \newline
\textbf{2. Complex and Kaehlerian Manifolds}\newline
Let $M$ configuration manifold$.$ A tensor field $J$ on $TM$ is called an
\textit{almost complex structure} on $TM$ if at every point $p$ of $TM,$ $J$
is endomorphism of the tangent space $T_{p}(TM)$ such that $\mathbf{J}^{2}=-%
\mathbf{I}.$ A manifold $TM$ with fixed almost complex structure $J$ is
called \textit{almost complex manifold}. If $(x^{i})$ and $(x^{i},\,y^{i})$
are coordinate systems of $M$ and $TM,$ then $\{\frac{\partial }{\partial
x^{i}},\frac{\partial }{\partial y^{i}}\}$ and $\{dx^{i},dy^{i}\}$ are
natural bases over $\mathbf{R}$ of the tangent space $T_{p}(TM)$ and the
cotangent space $T_{p}^{*}(TM)$ of $TM,$ respectively. Thus we get
\begin{equation}
\begin{array}{l}
J(\frac{\partial }{\partial x^{i}})=\frac{\partial }{\partial y^{i}},J(\frac{%
\partial }{\partial y^{i}})=-\frac{\partial }{\partial x^{i}}.
\end{array}
\label{2.1}
\end{equation}
Let $z^{i}=x^{i}+$\textbf{i}$\,y^{i},$ \textbf{i}$=\sqrt{-1},$ be a complex
local coordinate system of $TM.$ We define
\begin{equation}
\begin{array}{l}
\frac{\partial }{\partial z^{i}}=\frac{1}{2}\{\frac{\partial }{\partial x^{i}%
}-\mathbf{i}\frac{\partial }{\partial y^{i}}\},\,\,\frac{\partial }{\partial
\overline{z}^{i}}=\frac{1}{2}\{\frac{\partial }{\partial x^{i}}+\mathbf{i}%
\frac{\partial }{\partial y^{i}}\},\,\,dz^{i}=dx^{i}+\mathbf{i}dy^{i},\,\,d%
\overline{z}^{i}=dx^{i}-\mathbf{i}dy^{i},
\end{array}
\label{2.2}
\end{equation}
$\,\,\,\,\,\,\,$ $\,\,\,\,\,\,\,\,\,\,\,\,\,\,\,\,\,\,\,\,\,\,\,\,\,\,$
where $\frac{\partial }{\partial z^{i}}$ and $dz_{{}}^{i}$represent bases of
the tangent space $T_{p}(TM)$ and cotangent space $T_{p}^{*}(TM)$ of $TM$,
respectively. Then we calculate
\begin{equation}
\begin{array}{l}
J(\frac{\partial }{\partial z^{i}})=\mathbf{i}\frac{\partial }{\partial z^{i}%
},J(\frac{\partial }{\partial \overline{z}^{i}})=-\mathbf{i}\frac{\partial }{%
\partial \overline{z}^{i}}.
\end{array}
\label{2.3}
\end{equation}
\textit{Hermitian metric} on an almost complex manifold with almost complex
structure $J$ is a Riemannian metric $g$ on $TM$ such that
\begin{equation}
\begin{array}{l}
g(JX,Y)+g(X,JY)=0,\,\,\,\,\forall X,Y\in \chi (TM).
\end{array}
\label{2.4}
\end{equation}
An almost complex manifold $TM$ with a Hermitian metric is called an \textit{%
almost Hermitian manifold}. If $TM$ is a complex manifold, then $TM$ is
called a \textit{Hermitian manifold}. Let further $TM$ be a 2m-dimensional
almost Hermitian manifold with almost complex structure $J$ and Hermitian
metric $g$. The triple $(TM,J,g)$ is called an \textit{almost Hermitian
structure}. Let $(TM,J,g)$ be an almost Hermitian structure. The 2-form
defined by
\begin{equation}
\begin{array}{l}
\Phi (X,Y)=g(X,JY),\,\,\,\,\forall X,Y\in \chi (TM)
\end{array}
\label{2.5}
\end{equation}
is called\textit{\ the Kaehlerian form} of $(TM,J,g).$ An almost Hermitian
manifold is called \textit{almost Kaehlerian} if its Kaehlerian form\textbf{%
\ }$\Phi $ is closed. If, moreover, $TM$ is Hermitian, then $TM$ is called a
Kaehlerian manifold. \newline
\textbf{3. Complex Euler-Lagrange Equations}\newline
In this section, we deduce complex Euler-Lagrange equations for Classical
Mechanics structured on Kaehlerian manifold. Let $J$ be an almost complex
structure on the Kaehlerian manifold and $(z^{i},\overline{z}^{i})$ its
complex coordinates. The semispray $\xi $ and Liouville vector field $V=J\xi
$ on the Kaehlerian manifold are given by
\begin{equation}
\begin{array}{l}
\xi =\xi ^{i}\frac{\partial }{\partial z^{i}}+\overline{\xi }^{i}\frac{%
\partial }{\partial \overline{z}^{i}},\,\,J\xi =\mathbf{i}\xi ^{i}\frac{%
\partial }{\partial z^{i}}-\mathbf{i}\overline{\xi }^{i}\frac{\partial }{%
\partial \overline{z}^{i}}.
\end{array}
\label{3.1}
\end{equation}
We call \textit{the kinetic energy} and \textit{the potential energy of
system} the maps given by $T,P:TM\rightarrow \mathbf{C.}$ Then \textit{%
Lagrangian energy function} $L$ is the map $L:TM\rightarrow \mathbf{C}$ such
that $L=T-P$ and \textit{the energy function }$E_{L}$ associated $L$ is the
function given by $E_{L}=V(L)-L.$ The closed Kaehlerian form $\Phi _{L}$ is
the closed 2-form given by $\Phi _{L}=-dd_{J}L$ such that $d_{J}=\mathbf{i}%
\frac{\partial }{\partial z^{i}}dz^{i}-\mathbf{i}\frac{\partial }{\partial
\overline{z}^{i}}d\overline{z}^{i}:\mathcal{F}(TM)\rightarrow \wedge ^{1}TM.$
Then we have
\begin{equation}
\begin{array}{ll}
\Phi _{L}= & \mathbf{i}\frac{\partial ^{2}L}{\partial z^{j}\partial z^{i}}%
dz^{i}\wedge dz^{j}+\mathbf{i}\frac{\partial ^{2}L}{\partial \overline{z}%
^{j}\partial z^{i}}dz^{i}\wedge dz^{j}+\mathbf{i}\frac{\partial ^{2}L}{%
\partial z^{j}\partial \overline{z}^{i}}dz^{j}\wedge d\overline{z}^{i}+%
\mathbf{i}\frac{\partial ^{2}L}{\partial \overline{z}^{j}\partial \overline{z%
}^{i}}d\overline{z}^{j}\wedge d\overline{z}^{i}.
\end{array}
\label{3.2}
\end{equation}
Since the map $TM_{\Phi _{L}}:\chi (TM)\rightarrow \wedge ^{1}(TM)$ such
that $TM_{\Phi _{L}}(\xi )=i_{\xi }\Phi _{L}$ is an isomorphism, there
exists an unique vector $\xi $ on $TM$ such that the vector field $\xi $
holds the equality given by (\ref{1.1})$.$ Thus vector field $\xi $ on $TM$
is seen as a \textit{Lagrangian vector field} associated energy $L$\ on%
\textbf{\ }Kaehlerian manifold $TM$. Then
\begin{equation}
\begin{array}{ll}
i_{\xi }\Phi _{L}= & \mathbf{i}\xi ^{i}\frac{\partial ^{2}L}{\partial
z^{j}\partial z^{i}}dz^{j}-\mathbf{i}\xi ^{i}\frac{\partial ^{2}L}{\partial
z^{j}\partial z^{i}}\delta _{i}^{j}dz^{i}+\mathbf{i}\xi ^{i}\frac{\partial
^{2}L}{\partial \overline{z}^{j}\partial z^{i}}d\overline{z}^{j}-\mathbf{i}%
\overline{\xi }^{i}\frac{\partial ^{2}L}{\partial \overline{z}^{j}\partial
z^{i}}\delta _{i}^{j}dz^{i} \\
& +\mathbf{i}\xi ^{i}\frac{\partial ^{2}L}{\partial z^{j}\partial \overline{z%
}^{i}}d\overline{z}^{i}-\mathbf{i}\overline{\xi }^{i}\frac{\partial ^{2}L}{%
\partial z^{j}\partial \overline{z}^{i}}\delta _{i}^{j}dz^{j}+\mathbf{i}%
\overline{\xi }^{i}\frac{\partial ^{2}L}{\partial \overline{z}^{j}\partial
\overline{z}^{i}}\delta _{i}^{j}d\overline{z}^{i}-\mathbf{i}\overline{\xi }%
^{i}\frac{\partial ^{2}L}{\partial \overline{z}^{j}\partial \overline{z}^{i}}%
d\overline{z}^{j}.
\end{array}
\label{3.3}
\end{equation}
Since the closed Kaehlerian form $\Phi _{L}$ on $TM$ is symplectic
structure, we have
\begin{equation}
\begin{array}{l}
E_{L}=\mathbf{i}\xi ^{i}\frac{\partial L}{\partial z^{i}}-\mathbf{i}%
\overline{\xi }^{i}\frac{\partial L}{\partial \overline{z}^{i}}-L,
\end{array}
\label{3.4}
\end{equation}
and hence
\begin{equation}
\begin{array}{ll}
dE_{L}= & \mathbf{i}\xi ^{i}\frac{\partial ^{2}L}{\partial z^{j}\partial
z^{i}}dz^{j}-\mathbf{i}\overline{\xi }^{i}\frac{\partial ^{2}L}{\partial
z^{j}\partial \overline{z}^{i}}dz^{j}-\frac{\partial L}{\partial z^{j}}%
dz^{j}+\mathbf{i}\xi ^{i}\frac{\partial ^{2}L}{\partial \overline{z}%
^{j}\partial z^{i}}d\overline{z}^{j}-\mathbf{i}\overline{\xi }^{i}\frac{%
\partial ^{2}L}{\partial \overline{z}^{j}\partial \overline{z}^{i}}d%
\overline{z}^{j}-\frac{\partial L}{\partial \overline{z}^{j}}d\overline{z}%
^{j}.
\end{array}
\label{3.5}
\end{equation}
Considering \textbf{Eq.}(\ref{1.1}) and the integral curve $\alpha :\mathbf{C%
}\rightarrow TM$ of $\xi ,$i.e. $\xi (\alpha (t))=\frac{d\alpha (t)}{dt}$,
hence it is satisfied equations
\begin{equation}
\begin{array}{l}
-\mathbf{i}\left[ \xi ^{j}\frac{\partial ^{2}L}{\partial z^{j}\partial z^{i}}%
+\overline{\xi }^{i}\frac{\partial ^{2}L}{\partial \overline{z}^{j}\partial
z^{i}}\right] dz^{j}+\frac{\partial L}{\partial z^{j}}dz^{j}+\mathbf{i}%
\left[ \xi ^{j}\frac{\partial ^{2}L}{\partial z^{j}\partial \stackrel{.}{z}%
^{j}}+\overline{\xi }^{i}\frac{\partial ^{2}L}{\partial \overline{z}%
^{j}\partial \stackrel{.}{z}^{j}}\right] d\overline{z}^{j}+\frac{\partial L}{%
\partial \stackrel{.}{z}^{j}}d\overline{z}^{j}=0,
\end{array}
\label{3.6}
\end{equation}
where the dots mean derivatives with respect to the time. Then we have
\begin{equation}
\begin{array}{l}
\mathbf{i}\frac{d}{dt}\left( \frac{\partial L}{\partial z^{i}}\right) -\frac{%
\partial L}{\partial z^{i}}=0,\,\,\mathbf{i}\frac{d}{dt}\left( \frac{%
\partial L}{\partial \stackrel{.}{z}^{i}}\right) +\frac{\partial L}{\partial
\stackrel{.}{z}^{i}}=0,
\end{array}
\label{3.7}
\end{equation}
These equations infer \textit{complex Euler-Lagrange equations} whose
solutions are the paths of the semispray $\xi $ on Kaehlerian manifold $TM.$
Then $(TM,\Phi _{L},\xi )$ is a \textit{complex} \textit{Lagrangian system }%
on\textbf{\ }Kaehlerian manifold $TM.$\newline \textbf{Application
1:} \cite{Cabar}Let us consider the system illustrated in
\textbf{Figure1}. It consists of a light rigid rod of length $\ell
$, carrying a mass $m$ at one end, and hinged at the other end to
a vertical axis, so that it can swing freely in a vertical plane
and accelerate along the vertical axis. Let us obtain the
equations of motion for small oscillations by writing complex
Lagrange function. Complex Lagrangian function of the system is
\[
\begin{array}{l}
L=\frac{1}{2}m(\frac{\ell ^{2}(\stackrel{.}{z}+\stackrel{.}{\overline{z}}%
)^{2}}{A}+\frac{B^{2}}{4A}-\frac{1}{4}(\stackrel{.}{z}-\stackrel{.}{%
\overline{z}})^{2}+\frac{\ell \sin \theta (\stackrel{.}{z}+\stackrel{.}{%
\overline{z}})[\mathbf{i}(\stackrel{.}{z}-\stackrel{.}{\overline{z}})-\frac{B%
}{A^{1/2}}]}{A^{1/2}})-\frac{1}{2}\mathbf{i}mg(z-\overline{z})
\end{array}
\]
where $A=4\ell ^{2}-(z+\overline{z})^{2},\,B=(z+\overline{z})(\stackrel{.}{z}%
+\stackrel{.}{\overline{z}}).$ Then, considering \textbf{Eq.}(\ref{3.7}),
the complex-Lagrangian equations of the motion on the mechanical system, can
be calculated by
\[
\begin{array}{l}
\mathbf{i}\frac{d}{dt}Q-Q=0,\,\,\,\mathbf{i}\frac{d}{dt}W+W=0,
\end{array}
\]
such that
\[
\begin{array}{l}
Q=\frac{\ell ^{2}(\stackrel{.}{z}+\stackrel{.}{\overline{z}})B}{A^{2}}+\frac{%
(\stackrel{.}{z}+\stackrel{.}{\overline{z}})B}{4A}+\frac{(z+\overline{z}%
)B^{2}}{4A^{2}}+\frac{\ell B\sin \theta )[\mathbf{i}(\stackrel{.}{z}-%
\stackrel{.}{\overline{z}})-\frac{B}{A^{1/2}}]}{2A^{3/2}} \\
-\frac{\ell B\sin \theta (\stackrel{.}{z}+\stackrel{.}{\overline{z}})[%
\stackrel{.}{z}+\stackrel{.}{\overline{z}}+\frac{(z+\overline{z})B}{A}]}{2A}-%
\frac{1}{2}\mathbf{i}g
\end{array}
\]
and
\[
\begin{array}{l}
W=\frac{\ell ^{2}(\stackrel{.}{z}+\stackrel{.}{\overline{z}})}{A}+\frac{(%
\stackrel{.}{z}+\stackrel{.}{\overline{z}})B}{4A}-\frac{1}{4}\stackrel{.}{z}+%
\frac{1}{4}\stackrel{.}{\overline{z}}+\frac{\ell \sin \theta [\mathbf{i}(%
\stackrel{.}{z}-\stackrel{.}{\overline{z}})-\frac{B}{A^{1/2}}]}{2A^{1/2}}-%
\frac{\ell B\sin \theta (\stackrel{.}{z}+\stackrel{.}{\overline{z}})[\mathbf{%
i}-\frac{(z+\overline{z})}{A^{1/2}}]}{2A^{1/2}}.
\end{array}
\]
\textbf{Application 2: }\cite{Tekkoyun1}Let us consider the system
illustrated in \textbf{Figure2}. A central force field $f(\rho
)=A\rho ^{\alpha -1}(\alpha \neq 0,1)$ acts on a body with mass
$m$ in a constant gravitational field. Then let us find out the
Euler-Lagrange equations of the motion by assuming the body always
on the vertical plane.

The Lagrangian function of the system is,

\[
L=\frac{1}{2}m\stackrel{.}{z}\stackrel{.}{\overline{z}}-\frac{A}{\alpha }(%
\sqrt{z\overline{z}})^{\alpha }-\mathbf{j}mg\frac{(z-\overline{z})\sqrt{z%
\overline{z}}}{(z+\overline{z})\sqrt{1-\frac{(z-\overline{z})^{2}}{(z+%
\overline{z})^{2}}}}.
\]

Then, using \textbf{Eq.}(\ref{3.7}), the so-called Euler-Lagrange equations
of the motion on the mechanical systems, can be obtained, as follows:

\[
L1:\,\,\,\,\,\mathbf{i}\frac{\partial }{\partial t}S-S=0,\,\,\,\,\,\,\,\,\,%
\,\,L2:\,\,\,\,\,\,\,\mathbf{i}\frac{\partial }{\partial t}U+U=0,
\]

such that

\begin{eqnarray*}
S &=&-\frac{A}{2z}(\sqrt{z\overline{z}})^{\alpha }+\mathbf{i}\frac{mg(z-%
\overline{z})\overline{z}}{2\sqrt{z\overline{z}}(z+\overline{z})W}+\mathbf{i}%
\frac{mg\sqrt{z\overline{z}}}{(z+\overline{z})W} \\
&&-\mathbf{i}\frac{mg\sqrt{z\overline{z}}(z-\overline{z})}{(z+\overline{z}%
)^{2}W}-\mathbf{i}\frac{mg\sqrt{z\overline{z}}(z-\overline{z})(-\frac{(z-%
\overline{z})}{(z+\overline{z})^{2}}+\frac{(z-\overline{z})^{2}}{(z+%
\overline{z})^{3}})}{(z+\overline{z})W^{3}},
\end{eqnarray*}
and
\[
U=\frac{1}{2}m\stackrel{.}{\overline{z}}
\]

where $W=\sqrt{1-\frac{(z-\overline{z})^{2}}{(z+\overline{z})^{2}}}.$

\textbf{Conclusion}: In this study, the\textbf{\ }Lagrangian formalisms and
systems in Classical Mechanics had been intrinsically obtained making two
complex applications.


\begin{thebibliography}{9}
\bibitem{Deleon}  De Leon M., Rodrigues P.R., Generalized Classical
Mechanics and Field Theory, North-Holland Math. Stud., 112, Elsevier Sci.
Pub. Comp., Inc., New York, 1985.

\bibitem{Norbury}  Norbury, J.W., Lagrangians and Hamiltonians for High
School Students, arXiv: Physics/0004029v1, 2000.

\bibitem{Cabar}  Cabar G., Hamiltonian systems, Master Thesis, Pamukkale
University,Turkey, 2006.

\bibitem{Gibbons}  Gibbons G. W., Part III: Applications of Differential
Geometry to Physics, Cambridge CB3 0WA, UK., 2006.

\bibitem{Tekkoyun}  Tekkoyun M., G\"{o}rg\"{u}l\"{u} A., Higher Order
Complex Lagrangian and Hamiltonian Mechanics Systems'' , Physics Letters A,
vol.357, 261-269, 2006.

\bibitem{Tekkoyun1}  Tekkoyun M., Cabar G., Complex Lagrangians and
Hamiltonians, Journal of Arts and Sciences, \c{C}ankaya \"{U}niv.,
Fen-Ed.Fak., Issue 8/December 2007 .

\bibitem{Celik}  Celik B., Maple and Mathematics with Maple, Nobel
publication and distribution, Turkey, 2004.
\end{thebibliography}
\end{document}